\documentclass[a4paper]{amsart}
\usepackage{amsmath,amsfonts,amssymb}
\usepackage{amscd}
\usepackage{amsthm}
\usepackage{latexsym}
\input xy
\xyoption{all}

 \newtheorem{thm}{Theorem}[section]
 \newtheorem{cor}[thm]{Corollary}
 \newtheorem{lem}[thm]{Lemma}
 \newtheorem{prop}[thm]{Proposition}
 \theoremstyle{definition}
 
 \theoremstyle{remark}
 
 \numberwithin{equation}{section} 
 \newtheorem{ex}[thm]{Example}
 
\newlength{\defbaselineskip}
\newcommand{\setlinespacing}[1]%
           {\setlength{\baselineskip}{#1 \defbaselineskip}}

\begin{document}

\title []
 {To specify surfaces of revolution with pointwise 1-type Gauss map
in 3-dimensional Minkowski space}
\author [Milani] {V. Milani}
\author [Shojaei-Fard] {A. Shojaei-Fard}

\address{Department of Mathematics, Shahid Beheshti University, 1983963113 Tehran, Iran}
\email{v-milani@cc.sbu.ac.ir} 
\email{a\_shojaei@sbu.ac.ir}


\date{}

\dedicatory{}

\commby{}

%%% ----------------------------------------------------------------------

%%% ----------------------------------------------------------------------
\maketitle 
\setlinespacing{1.12}
%%% ----------------------------------------------------------------------

\begin{abstract}
In this paper, by the studying of the Gauss map, Laplacian operator, curvatures of surfaces in $\mathbb{R}_{1}^{3}$ and Bour's theorem, we are going to identify surfaces of revolution with pointwise 1-type Gauss map property in $3-$dimensional Minkowski space.
\end{abstract}

\section*{Introduction}

The classification of submanifolds in Euclidean and Non-Euclidean
spaces is one of the interesting topics in differential geometry and in this way one can find some
attempts in terms of finite type submanifolds \cite{BCV1, C1, C2, C3,
CP1}. On the other hand  Kobayashi in \cite{K1} classified space-like
ruled minimal surfaces in $\mathbb{R}_{1}^{3}$ and its extension to
the Lorentz version is done by de Woestijne in \cite{W1}. In continue, people encounter with the following 
problem:

\emph{Classify all surfaces in 3-dimensional Minkowski space
satisfying the pointwise 1-type Gauss map condition $\Delta N = kN$
for the Gauss map $N$ and some function $k$.}

In 2000, D.W.Yoon and Y.H.Kim in \cite{KY1} classified minimal ruled
surfaces in terms of pointwise 1-type Gauss map in
$\mathbb{R}_{1}^{3}$.

On suitability oriented surface $M$ in $\mathbb{R}^{3}$ with
positive Gaussian curvature $K$, one can induce a positive definite
second fundamental form $II$ with component functions $e$, $f$, $g$.
The second Gaussian curvature is defined by

\begin{equation}
K_{II}= \frac{1}{(|eg|-f^{2})^{2}}( \left|
         \begin{array}{rrr}
              -\frac{1}{2}e_{tt}+f_{st}-\frac{1}{2}g_{ss} & \frac{1}{2}e_{s} & f_{s}-\frac{1}{2}e_{t}
              \\
              f_{t}-\frac{1}{2}g_{s} &  e &  f\\
              \frac{1}{2}g_{t} &  f  &  g
          \end{array} \right| -\left|
         \begin{array}{rrr}
              0 & \frac{1}{2}e_{t} & \frac{1}{2}g_{s}
              \\
              \frac{1}{2}e_{t} &  e &  f\\
              \frac{1}{2}g_{s} &  f  &  g
          \end{array} \right| ).\ 
\end{equation} \cite{KY2}

We can extend it to the surfaces in
$\mathbb{R}_{1}^{3}$.

In 2004, D.W.Yoon and Y.H.Kim in  \cite{KY2},
classified ruled surfaces in terms of the second Gaussian curvature,
the mean curvature and the Gaussian curvature in 3-dimensional
Minkowski space. On the other hand, in 2001, Ikawa in \cite{I1}
proved Bour's theorem in
$\mathbb{R}_{1}^{3}$. He showed that
 
\begin{quote}
\emph{A generalized helicoid is isometric to a surface of revolution
in $\mathbb{R}_{1}^{3}$.}
\end{quote}

In this paper, the above problem is answered for the surfaces of
revolution in 3-dimensional Minkowski space.

\section{Classification}

Let $\mathbb{R}_{1}^{3}$ be a 3-dimensional Minkowski space with the
scalar product $\langle \ , \ \rangle$ of index 1 defined as
\begin{equation}
\langle X,Y\rangle = X_{1}Y_{1} +X_{2}Y_{2} - X_{3}Y_{3}
\end{equation}

for every vectors $X=(X_{i})$ and $Y=(Y_{i})$ in
$\mathbb{R}_{1}^{3}$. 

A vector $X$ of $\mathbb{R}_{1}^{3}$ is said
to be {\it space-like} if $\langle X,X \rangle >0$ or $X=0$, {\it time-like} if
$\langle X,X \rangle <0$ and {\it light-like or null} if $\langle X,X
\rangle =0$ and $X\neq 0$. A time-like or light-like vector in
$\mathbb{R}_{1}^{3}$ is said to be {\it causal}.

\begin{lem} \label{1}
There are no causal vectors in
$\mathbb{R}_{1}^{3}$ orthogonal to a time-like vector \cite{G1}.
\end{lem}

A lorentz cross product $X\times Y$ is given by

\begin{equation}
X\times Y=(X_{2}Y_{3}-X_{3}Y_{2} , X_{3}Y_{1}-X_{1}Y_{3},X_{2}Y_{1}-X_{1}Y_{2} ).
\end{equation}

A curve in $\mathbb{R}_{1}^{3}$ is called space-like, time-like or
light-like if the tangent vector at any point is space-like,
time-like or light-like, respectively. A plane in
$\mathbb{R}_{1}^{3}$ is space-like, time-like or light-like if its
Euclidean unit normal is time-like, space-like or light-like,
respectively. A surface in $\mathbb{R}_{1}^{3}$ is space-like,
time-like or light-like if the tangent plane at any point is
space-like, time-like or light-like, respectively. Let $M$ be a
surface in $\mathbb{R}_{1}^{3}$. The Gauss map $N:M\longrightarrow
Q^{2}(\epsilon) \subset \mathbb{R}_{1}^{3}$ which sends each point
of $M$ to the unit normal vector to $M$ at that point is called the {\it
Gauss map of surface $M$}. Here $\epsilon(=\pm1)$ denotes the sign of
the vector field $N$ and $Q^{2}(\epsilon)$ is a 2-dimensional space
form as follows:

\begin{equation}
Q^{2}(\epsilon) =
\end{equation}
$$
 \left\lbrace
  \begin{array}{c l}
    S^{2}_{1}(1) & \text{in $\mathbb{R}_{1}^{3} (\epsilon = 1)$} \\
    H^{2}(-1) &  \text{in $\mathbb{R}_{1}^{3} (\epsilon = -1)$}
  \end{array}
\right. 
$$

It is well known that in terms of local coordinates $(x_{i})$ of $M$
the Laplacian can be written by

\begin{equation} \label{*}
\Delta =-\frac{1}{\sqrt{|\textbf{G}|}} \sum_{i,j} \frac{\partial}{\partial
x^{i}}(\sqrt{|\textbf{G}|}g^{ij}\frac{\partial}{\partial x^{j}})
\end{equation}

where $\textbf{G}=\det(g_{ij})$, $(g^{ij})=(g_{ij})^{-1}$ and
$(g_{ij})$ are the components of the metric of $M$ with respect to
$(x_{i})$.

Now, we define a ruled surface $M$ in a three-dimensional Minkowski
space $\mathbb{R}_{1}^{3}$. Let $I$ be an open interval in the real
line $\mathbb{R}$. Let $\alpha=\alpha(s)$ be a curve in
$\mathbb{R}_{1}^{3}$ defined on $I$ and $\beta=\beta(s)$ a
transversal vector field along $\alpha$. For an open interval $J$ in
$\mathbb{R}$, let $M$ be a ruled surface parameterized by:

\begin{equation}
x=x(s,t)=\alpha(s)+t\beta(s)\ \ s\in I,\ t\in J.
\end{equation}

According to the character of the base curve $\alpha$ and the
director curve $\beta$ the ruled surfaces are classified into the
following six groups.

If the base curve $\alpha$ is space-like or time-like, then the
ruled surface $M$ is said to be of { \it type $M_{+}$ or type $M_{-}$},
respectively. Also the ruled surface of type $M_{+}$ can be divided
into three types. When $\beta$ is space-like, it is said to be of
{\it type $M^{1}_{+}$} or {\it $M^{2}_{+}$} if $\beta'$ is non-null or
light-like, respectively. By \ref{1} when $\beta$ is time-like, $\beta'$ must
be space-like. In this case, $M$ said to be
of {\it type $M^{3}_{+}$}. On the other hand, for the ruled surface of
type $M_{-}$, it is also said to be of {\it type $M^{1}_{-}$} or
{\it  $M^{2}_{-}$} if $\beta'$ is non-null or light-like, respectively.
Note that in the case of type $M_{-}$ the director curve $\beta$ is
always space-like. The ruled surface of type $M^{1}_{+}$ or
$M^{2}_{+}$ (resp. $M^{3}_{+}$, $M^{1}_{-}$, $M^{2}_{-}$) is clearly
space-like (resp. time-like). If the base curve $\alpha$ and the
director curve $\beta$ are light-like, then the ruled surface is
called {\it null scroll} and it is a time-like surface.

Now we modify the definition of a surface of revolution and generalized
helicoid in $\mathbb{R}_{1}^{3}$. For an open interval $I\subset
\mathbb{R}$, let $\gamma:I\longrightarrow \Pi$ be a curve in a plane
$\Pi$ in $\mathbb{R}_{1}^{3}$ ({\it profile curve}) and let $l$ be a
straight line in $\Pi$ which does not intersect the curve $\gamma$
({\it axis}). A surface of revolution in $\mathbb{R}_{1}^{3}$ is
defined by the Lorentzian rotation $\gamma$ around $l$.

Suppose the case when a profile curve $\gamma$ rotates around the
axis $l$, it simultaneously displaces parallel to $l$ so that the
speed of displacement is proportional to the speed of rotation. Then
the resulting surface is called {\it generalized helicoid}.

In the following at first we give some examples of surfaces of revolution in
$\mathbb{R}_{1}^{3}$ and next will use them.

\begin{ex} \label{2}
For the constants $a,b$, let $M$ be a surface in $\mathbb{R}_{1}^{3}$ with the parametric representation 
\begin{equation}
R(s,t)=( \sqrt{(t+a)^{2}-b^{2}}\cos s , \sqrt{(t+a)^{2}-b^{2}}\sin
s , \int\sqrt{\frac{b^{2}}{(t+a)^{2}-b^{2}}}dt  )
\end{equation}
where
$b^{2}<(t+a)^{2}$. It is called {\it surface of
revolution of the 1st kind as space-like}.
\end{ex}

\begin{ex}  \label{3}
For the constants $a,b$, let $M$ be a surface in $\mathbb{R}_{1}^{3}$ with the parametric representation 
\begin{equation}
R(s,t)=(\sqrt{b^{2}-(t+a)^{2}}\sinh s ,
-b\sin^{-1}(\frac{\sqrt{b^{2}-(t+a)^{2}}}{-b}) ,
\sqrt{b^{2}-(t+a)^{2}}\cosh s)
\end{equation}
where $(t+a)^{2}<b^{2}$. It is called {\it surface of revolution of the 2nd kind as
space-like}.
\end{ex}

\begin{ex} \label{4}
For constants $a,b$, let $M$ be a surface in $\mathbb{R}_{1}^{3}$ with the parametric representation 
\begin{equation}
R(s,t)=( \sqrt{(t+a)^{2}-b^{2}}\cosh s ,-b\cosh^{-1}(\frac{\sqrt{(t+a)^{2}-b^{2}}}{-b}) ,\sqrt{(t+a)^{2}-b^{2}}\sinh
 s)
\end{equation}
where $b^{2}<(t+a)^{2}$. It is called {\it surface of
revolution of the 2nd kind as time-like}.
\end{ex}

\begin{ex} \label{5}
For constants $a,b$, let $M$ be a surface in $\mathbb{R}_{1}^{3}$ with the parametric representation 
\begin{equation}
R(s,t)=(\int \sqrt{\frac{-b^{2}} { b^{2}+(t+a)^{2}} } dt , \sqrt{b^{2}+(t+a)^{2}}\sinh
s , \sqrt{b^{2}+(t+a)^{2}}\cosh s).
\end{equation}
It is called
{\it surface of revolution of the 3rd kind as Lorentzian}.
\end{ex}

\begin{prop} \label{6}
Let $M$ be a helicoid of the 1st kind as
space-like
$$x(s,t)=((t+a)\cos s , (t+a)\sin s , -bs) $$
where $|a|>|b|>0$, $t<min(-a-b,-a+b)$ or $t>max(-a-b,-a+b)$. This
surface is isometric to a minimal surface of revolution with
pointwise 1-type Gauss map property.
\end{prop}

\begin{proof}

According to the Bour's theorem in Minkowski 3-space \cite{I1},
for each helicoidal surface there exists an isometric surface of revolution to it. Therefore one can
see that helicoid of the 1st kind as space-like is isometric to the
surface of revolution of the 1st kind as space-like. By the parametrization of this surface of revolution, its Gauss map is given by

$$N=\frac{R_{s}\times R_{t}}{\parallel R_{s}\times R_{t} \parallel
}=-\frac{1}{\sqrt{(t+a)^{2}-b^{2}}}(b\cos s , b\sin s ,t+a). $$ The
components $(g_{ij})$ of the metric with respect to the first
fundamental forms of this surface are

$$E=g_{11}=\langle R_{s},R_{s}\rangle=(t+a)^{2}-b^{2},$$
$$F=g_{12}=\langle R_{s},R_{t}\rangle=0, $$
$$F=g_{21}= \langle R_{t},R_{s}\rangle=0, $$
$$G=g_{22}=\langle R_{t},R_{t}\rangle=1.$$
By (\ref{*}),
$$\Delta N=2b^{2}( (t+a)^{2}-b^{2} ) ^{-\frac{5}{2}}(b\cos s , b\sin s ,t+a). $$
Then for some function $k$, $\Delta N = kN$ such that $k=-2b^{2}(
(t+a)^{2}-b^{2} ) ^{-2}$. In other words, this surface of revolution
has pointwise 1-type Gauss map property. On the other hand, the
second fundamental forms of surface of revolution of the 1st kind as
space-like are

$$e=\langle R_{ss},N\rangle=b, $$
$$f=\langle R_{st},N\rangle=\langle R_{ts},N\rangle,$$
$$g=\langle R_{tt},N\rangle=-b( (t+a)^{2}-b^{2} )^{-1}.$$
The mean curvature $H$ is given by
$$ H=\frac{Eg-2Ff+Ge}{2|EG-F^{2}|}=0. $$
Therefore surface of revolution of the 1st kind as space-like is a
maximal surface and its Gauss map is of pointwise 1-type.

\end{proof}

\begin{prop} \label{7}
Let $M$ be a helicoid of the 2nd kind as
space-like
$$x(s,t)=((t+a)\cosh s , -bs , (t+a)\sinh s ) $$
where $|b|>|a|$, $min(-a-b,-a+b)<t<max(-a-b,-a+b)$. This surface is
isometric to a minimal surface of revolution with pointwise 1-type
Gauss map property.
\end{prop}

\begin{proof}

According to the Bour's theorem in Minkowski 3-space \cite{I1},
for every helicoidal surface one can find its isometric surface of revolution. Therefore the helicoid of the 2nd kind as space-like is isometric to the surface of revolution of the 2nd kind as space-like. By the parametrization of this surface of revolution, its Gauss map is given by
$$N=\frac{R_{s}\times R_{t}}{\parallel R_{s}\times R_{t} \parallel
}=\frac{1}{\sqrt{b^{2}-(t+a)^{2}}}( -b\sinh s ,t+a , -b\cosh s).
$$
The components $(g_{ij})$ of the metric with respect to the first
fundamental forms of this surface are
$$E=g_{11}=\langle R_{s},R_{s}\rangle=b^{2}-(t+a)^{2},$$
$$F=g_{12}=\langle R_{s},R_{t}\rangle=0, $$
$$F=g_{21}= \langle R_{t},R_{s}\rangle=0, $$
$$G=g_{22}=\langle R_{t},R_{t}\rangle=1.$$
By (\ref{*}),
$$\Delta N =-2b^{2}(b^{2}-(t+a)^{2})^{-\frac{5}{2}}( -b\sinh s ,t+a , -b\cosh s).$$
Then for some function $k$, $\Delta N = kN$  such that
$k=-2b^{2}(b^{2}-(t+a)^{2})^{-2}$. In other words, this surface of
revolution has pointwise 1-type Gauss map property. And also, the second fundamental forms of surface of revolution of the 2nd kind as space-like are
$$e=\langle R_{ss},N\rangle=b, $$
$$f=\langle R_{st},N\rangle=\langle R_{ts},N\rangle=0,$$
$$g=\langle R_{tt},N\rangle=-b(b^{2}-(t+a)^{2})^{-1}.$$
The mean curvature $H$ is given by
$$ H=\frac{Eg-2Ff+Ge}{2|EG-F^{2}|}=0. $$
Therefore surface of revolution of the 2nd kind as space-like is a
maximal surface and its Gauss map has pointwise 1-type property.

\end{proof}

\begin{prop} \label{8}
Let $M$ be a helicoid of the 2nd kind as
time-like
$$x(s,t)=(  (t+a)\cosh s ,-bs , (t+a)\sinh s ) $$
where $|a|>|b|>0$, $t<min(-a-b,-a+b)$ or $ t>max(-a-b,-a+b)$. This
surface is isometric to a minimal surface of revolution with
pointwise 1-type Gauss map property.
\end{prop}

\begin{proof}

By Bour's theorem in Minkowski 3-space,
the helicoid of the 2nd kind as time-like is isometric
to the surface of revolution of the 2nd kind as time-like. The
Gauss map of this surface of revolution is
$$N=\frac{R_{s}\times R_{t}}{\parallel R_{s}\times R_{t} \parallel
}=\frac{1}{\sqrt{(t+a)^{2}-b^{2}}}( -b\sinh s ,t+a , -b\cosh s).
$$
The components $(g_{ij})$ of the metric with respect to the first
fundamental forms of this surface are
$$E=g_{11}=\langle R_{s},R_{s}\rangle=-((t+a)^{2}-b^{2}),$$
$$F=g_{12}=\langle R_{s},R_{t}\rangle=0, $$
$$F=g_{21}= \langle R_{t},R_{s}\rangle=0, $$
$$G=g_{22}=\langle R_{t},R_{t}\rangle$$
$$=(t+a)^{2}((t+a)^{2}-b^{2})^{-1}+(t+a)^{2}((t+a)^{2}-b^{2})^{-1}sinh^{2}(\frac{((t+a)^{2}-b^{2})^{\frac{1}{2}}}{-b})cosh^{-4}(\frac{((t+a)^{2}-b^{2})^{\frac{1}{2}}}{-b}).$$
By (\ref{*}),
$$\Delta N =-2b^{2}((t+a)^{2}-b^{2})^{-\frac{5}{2}}( -b\sinh s ,t+a , -b\cosh s).$$
Then for some function $k$, $\Delta N = kN$ such that
$k=-2b^{2}((t+a)^{2}-b^{2})^{-2}$. In other words, this surface of
revolution has pointwise 1-type Gauss map property. On the other
hand by the calculating of the second fundamental forms of the surface of
revolution of the 2nd kind as time-like (similar to the surface of
revolution of the 2nd kind as space-like) \cite{I1}, one can see
that the mean curvature $H$ is given by
$$ H=\frac{Eg-2Ff+Ge}{2|EG-F^{2}|}=0. $$
Therefore surface of revolution of the 2nd kind as time-like is a
minimal surface and its Gauss map is of pointwise 1-type.

\end{proof}

\begin{prop} \label{9}
Let $M$ be a helicoid of the 3rd kind as
Lorentzian
$$x(s,t)=( bs , (t+a)\sinh s , (t+a)\cosh s ) $$
where $|a|<|b|$, $min(-a-b,-a+b)<t<max(-a-b,-a+b)$. This surface is
isometric to a minimal surface of revolution with pointwise 1-type
Gauss map property.
\end{prop}

\begin{proof}

Bour's theorem gives an isometric between  the helicoid of the 3rd kind as Lorentzian and surface of revolution of the 3rd kind as Lorentzian. Its
Gauss map is given by 
$$N=\frac{R_{s}\times R_{t}}{\parallel R_{s}\times R_{t} \parallel
}=-((t+a)^{2}+b^{2})^{-\frac{1}{2}}(t+a , ib\sinh s , ib\cosh s).
$$
The components $(g_{ij})$ of the metric with respect to the first
fundamental forms of this surface are

$$E=g_{11}=\langle R_{s},R_{s}\rangle=(t+a)^{2}+b^{2},$$
$$F=g_{12}=\langle R_{s},R_{t}\rangle=0, $$
$$F=g_{21}= \langle R_{t},R_{s}\rangle=0, $$
$$G=g_{22}=\langle R_{t},R_{t}\rangle=-1.$$
By (\ref{*}),
$$\Delta N =2b^{2}((t+a)^{2}+b^{2})^{-\frac{5}{2}}(t+a , ib\sinh s , ib\cosh s).$$
Then $\Delta N = kN$ for some function $k$ such that
$k=-2b^{2}((t+a)^{2}+b^{2})^{-2}$. It means that, this surface of
revolution has pointwise 1-type Gauss map property. 

The second fundamental forms of the surface of revolution of the
3rd kind as Lorentzian are
$$e=\langle R_{ss},N\rangle=ib, $$
$$f=\langle R_{st},N\rangle=\langle R_{ts},G\rangle,$$
$$g=\langle R_{tt},N\rangle=ib((t+a)^{2}+b^{2})^{-1}.$$
The mean curvature $H$ is given by
$$ H=\frac{Eg-2Ff+Ge}{2|EG-F^{2}|}=0. $$
Therefore surface of revolution of the 3rd kind as Lorentzian is a
minimal surface and its Gauss map is of pointwise 1-type.
\end{proof}

\begin{prop} \label{10}
Let $M$ be a conjugate of Enneper's
surface of the 2nd kind as space-like or time-like
$$x(s,t)=(hs^{2}+t,h(\frac{s^{3}}{3}-s)+ts, h(\frac{s^{3}}{3}+s)+ts). $$
This surface is isometric to a minimal or maximal surface of
revolution with pointwise 1-type Gauss map property.
\end{prop}

\begin{proof}

According to the Bour's theorem in Minkowski 3-space
\cite{I1,W1}, we can see that the conjugate of Enneper's surface of
the 2nd kind as space-like or time-like is isometric to minimal or
maximal surface of revolution Enneper of the 2nd, 3rd kind
$$x(s,t)=(at^{3}+t-s^{2}t+b,-2st,at^{3}-t-s^{2}t+b) $$
or
$$x(s,t)=(-at^{3}+t-s^{2}t+b,-2st,-at^{3}-t-s^{2}t+b) $$
respectively, where $a>0,b \in \mathbb{R}$ and $t\neq 0$.

On the other hand, the conjugate of Enneper's surface of the 2nd kind as
space-like and the surface of revolution Enneper of the 2nd kind have
the same Gauss map. Also, conjugate of Enneper's surface of the 2nd
kind as time-like and surface of revolution Enneper of the 3rd kind
have the same Gauss map. Since the conjugate of Enneper's surface of the
2nd kind as space-like or time-like has pointwise 1-type Gauss map
property \cite{KY1}, the proof is completed.

\end{proof}

\begin{prop} \label{11}
 Let $M$ be a non-developable ruled surface of
type $M^{1}_{+}$ or $M^{3}_{+}$ in $\mathbb{R}_{1}^{3}$, such that
has one of the following conditions:
$$aK_{II}+bH=constant,\ a,b\in
\mathbb{R}-\{0\} ,\ 2a-b\neq 0,\ along\ each \ ruling$$ or
$$aH+bK=constant,\ a\neq 0,\ b\in \mathbb{R},\ along\ each \
ruling$$ or
$$aK_{II}+bK=constant,\ a\neq 0,\ b\in \mathbb{R},\
along\ each \ ruling.$$ Then $M$ is isometric to an open part of one
of the following surfaces of revolution:
\begin{enumerate}
  \item Surface of revolution of the 1st kind as space-like,
  \item Surface of revolution of the 2nd kind as space-like,
  \item Surface of revolution of the 3rd kind as Lorentzian.
\end{enumerate}
\end{prop}

\begin{proof}

According to theorems 4.1, 4.2 and 4.4 in \cite{KY2}, we know
that every non-developable ruled surfaces of type $M^{1}_{+}$ or
$M^{3}_{+}$ in $\mathbb{R}_{1}^{3}$ such that have one of the above
conditions, are open parts of one of the following surfaces:

\begin{quote}
   The helicoid of the 1st, 2nd kind as space-like,
\end{quote}

\begin{quote}
   The helicoid of the 3rd kind as Lorentzian.
\end{quote}
On the other hand \ref{6}, \ref{7} and \ref{9} show that
these surfaces are isometric to

\begin{quote}
  Surface of revolution of the 1st, 2nd kind as space-like,
\end{quote}

\begin{quote}
  Surface of revolution of the 3rd kind as Lorentzian,
\end{quote}
respectively.

\end{proof}

\begin{prop} \label{12}
Let $M$ be a non-developable ruled surface of
type $M^{1}_{-}$ and $M$ isn't an open part of the helicoid of the
1st kind as time-like in $\mathbb{R}_{1}^{3}$ such that this ruled
surface satisfies one of the following conditions:
$$aK_{II}+bH=constant,\ a,b\in \mathbb{R}-\{0\},\ 2a-b\neq 0,\
along\ each\ ruling$$ or
$$aH+bK=constant,\ a\neq 0,\ b\in
\mathbb{R},\ along\ each\ ruling$$ or
$$aK_{II}+bK=constant,\ a\neq
0,\ b\in \mathbb{R},\ along\ each\ ruling.$$ Then $M$ is isometric
to an open part of the following surface:
\begin{center}
Surface of revolution of the 2nd kind as time-like.
\end{center}
\end{prop}
\begin{proof}
According to theorems 4.1, 4.2 and 4.4 in \cite{KY2}, we know
that every non-developable ruled surfaces of type $M^{1}_{-}$ and
not an open part of the helicoid of the 1st kind as time-like in
$\mathbb{R}_{1}^{3}$ such that have one of the above conditions, are
open parts of the helicoid of the 2nd kind as time-like. On the
other hand \ref{8} shows that this surface is isometric to
the surface of revolution of the 2nd kind as time-like. This completes
the proof.
\end{proof}

\begin{prop} \label{13}
Let $M$ be a non-developable ruled surface of
type $M^{2}_{+}$ or $M^{2}_{-}$ in $\mathbb{R}_{1}^{3}$, such that
has one of the following conditions:
$$aK_{II}+bH=constant,\ a,b\in
\mathbb{R}-\{0\} ,\ 2a-b\neq 0,\ along\ each \ ruling$$ or
$$aH+bK=constant,\ a\neq 0,\ b\in \mathbb{R},\ along\ each \
ruling.$$ Then $M$ is isometric to an open part of the following
surfaces of revolution:
\begin{center}
  \item surfaces of  Enneper of the 2nd, 3rd kind.
\end{center}
\end{prop}

\begin{proof}

By theorems 4.1 and 4.2 in \cite{KY2}, one can see that
every non-developable ruled surfaces of type $M^{2}_{+}$ or
$M^{2}_{-}$ in $\mathbb{R}_{1}^{3}$ such that have one of the above
conditions, are open parts of the conjugate of Enneper's surfaces of
the 2nd kind as space-like or time-like. On the other hand \ref{10} shows that these surfaces isometric to surfaces
of Enneper of the 2nd, 3rd kind.

\end{proof}

\begin{cor} \label{14}
Let $R$ be a surface of revolution such
that is isometric to an open part of the non-developable ruled surface
$M$ of type $M^{1}_{+}$ or $M^{3}_{+}$ where $M$ satisfies one
of the following conditions:
$$aK_{II}+bH=constant,\ a,b\in
\mathbb{R}-\{0\},\ 2a-b\neq 0,\ along\ each\ ruling$$ or
$$aH+bK=constant,\ a\neq 0,\ b\in \mathbb{R},\ along\ each\
ruling$$ or
$$aK_{II}+bK=constant,\ a\neq 0,\ b\in \mathbb{R},\
along\ each\ ruling.$$ Then $R$ has pointwise 1-type Gauss map
property.
\end{cor}

\begin{proof}
According to \ref{11}, $R$ is an open part of one of the
following surfaces of revolution:

\begin{quote}
   Surface of revolution of the 1st kind as space-like,
\end{quote}

\begin{quote}
   Surface of revolution of the 2nd kind as space-like,
\end{quote}

\begin{quote}
   Surface of revolution of the 3rd kind as Lorentzian.
\end{quote}
On the other hand \ref{6}, \ref{7} and \ref{9} show that
these surfaces have pointwise 1-type Gauss map property.
\end{proof}

\begin{cor} \label{15}
Let $R$ be a surface of revolution such
that is isometric to an open part of the non-developable ruled surface
$M$ of type $M^{1}_{-}$ and $M$ isn't an open part of the helicoid
of the 1st kind as time-like in $\mathbb{R}_{1}^{3}$ and also this
ruled surface satisfies one of the following conditions:
$$aK_{II}+bH=constant,\ a,b\in \mathbb{R}-\{0\},\ 2a-b\neq 0,\ along
\ each\ ruling$$ or
$$aH+bK=constant,\ a\neq 0,\ b\in \mathbb{R},\
along\ each\ ruling$$ or
$$aK_{II}+bK=constant,\ a\neq 0,\ b\in \mathbb{R},\ along\ each\
ruling.$$ Then $R$ has pointwise 1-type Gauss map property.
\end{cor}

\begin{proof}
According to \ref{12}, $R$ is an open part of the surface of
revolution of the 2nd kind as time-like. On the other hand
\ref{8} shows that this surface has pointwise 1-type Gauss
map property.
\end{proof}

\begin{cor} \label{16}
Let $R$ be a surface of revolution such
that is isometric to an open part of the non-developable ruled surface
$M$ of type $M^{2}_{+}$ or $M^{2}_{-}$ where $M$ satisfies one
of the following conditions:
$$aK_{II}+bH=constant,\ a,b\in
\mathbb{R}-\{0\},\ 2a-b\neq 0,\ along\ each\ ruling$$ or
$$aH+bK=constant,\ a\neq 0,\ b\in \mathbb{R},\ along\ each\
ruling.$$  Then $R$ has pointwise 1-type Gauss map property.
\end{cor}

\begin{proof}
According to \ref{13}, $R$ is an open part of the surfaces
Enneper of the 2nd and 3rd kind. On the other hand \ref{10} shows that these surfaces have pointwise 1-type Gauss map property.
\end{proof}

\begin{cor} \label{17}
Let $R$ be one of the following surfaces:
\begin{enumerate}

\item Surface of revolution of the 1st kind as space-like,
\item Surface of revolution of the 2nd kind as space-like or
time-like,
\item Surface of revolution of the 3rd kind as Lorentzian,
\item Surfaces of Enneper of the 2nd, 3rd kind.

\end{enumerate}
Then $R$ is a part of one of the following surfaces:
\begin{enumerate}

  \item A space-like or time-like plane,
  \item The catenoids of the 1st, 2nd, 3rd, 4th, 5th kind,
  \item Surfaces of Enneper of the 2nd, 3rd kind.

\end{enumerate}
\end{cor}

\begin{proof}
By \ref{6}, \ref{7}, \ref{8}, \ref{9} and also \ref{10}, these
surfaces of revolution are minimal or maximal. On the other hand by
\cite{W1}, we know that the only minimal or maximal surfaces of
revolution in Minkowski 3-space are an open part of one of the
following surfaces: A space-like or time-like plane, The catenoids
of the 1st, 2nd, 3rd, 4th, 5th kind, Surfaces of Enneper of the 2nd,
3rd kind.
\end{proof}

And finally one can reach to a nice characterization of the minimal and maximal surfaces of revolution by the pointwise 1-type property. We have

\begin{prop} \label{18} 
Let $R$ be a surface of revolution in
$\mathbb{R}_{1}^{3}$. $R$ has pointwise 1-type Gauss map property if
and only if $R$ be an open part of the minimal or maximal surfaces
of revolution.
\end{prop}

% ------------------------------------------------------------------------
%Included for Gather Purpose only:
%input "Xbib.bib"
%\bibliographystyle{amsplain}
%\bibliography{xbib}
\end{document}